\documentclass[12pt]{amsart}

\usepackage[cp1251]{inputenc}
\usepackage{amsthm,amssymb,amsmath,amsfonts}
\usepackage{graphicx}
\usepackage[matrix,arrow,curve]{xy}
\newcommand{\CC}{\ensuremath{\mathbb{C}}}

\newcommand{\CG}{\ensuremath{\mathfrak{C}}}

\newcommand{\Q}{\ensuremath{\mathbb{Q}}}

\newcommand{\Z}{\ensuremath{\mathbb{Z}}}

\newcommand{\A}{\ensuremath{\mathbb{A}}}

\newcommand{\F}{\ensuremath{\mathbb{F}}}

\newcommand{\XX}{\ensuremath{\overline{X}}}

\newcommand{\ka}{\ensuremath{\Bbbk}}

\newcommand{\kka}{\ensuremath{\overline{\Bbbk}}}

\newcommand{\K}{\ensuremath{\mathrm{K}}}

\newcommand{\Pro}{\ensuremath{\mathbb{P}}}

\newcommand{\Aut}{\ensuremath{\operatorname{Aut}}}

\newcommand{\Pic}{\ensuremath{\operatorname{Pic}}}

\newcommand{\Char}{\ensuremath{\operatorname{char}}}

\newcommand{\ord}{\ensuremath{\operatorname{ord}}}

\newcommand{\matr}[9]{\ensuremath{\left(
\begin{array}{ccc} #1&#2&#3\\#4&#5&#6\\#7&#8&#9\\
\end{array}
\right)}}

\makeatletter
\@addtoreset{equation}{section}
\makeatother

\newtheorem{theorem}[equation]{Theorem}
\newtheorem{proposition}[equation]{Proposition}
\newtheorem{lemma}[equation]{Lemma}
\newtheorem{corollary}[equation]{Corollary}

\theoremstyle{definition}

\newtheorem{definition}[equation]{Definition}

\theoremstyle{remark}
\newtheorem{remark}[equation]{Remark}
\newtheorem{convention}[equation]{Convention}

\title{Rationality of the quotient of $\Pro^2$ by finite group of automorphisms over arbitrary field of characteristic zero}
\address{Department of Algebra, Faculty of Mathematics, MSU, Moscow 117234, Russia}
\address{Laboratory of Algebraic Geometry, GU-HSE, 7 Vavilova str., Moscow 117312, Russia}
\address{Independent University of Moscow, Moscow 119002, Russia}
\email{trepalin@mccme.ru}
\thanks{The author was partially supported by AG Laboratory HSE, RF government grant, ag. 11.G34.31.0023 and the grants RFFI 11-01-00336-a, MK-503.2010.1 and N.SH.-4713.2010.1}
\author{Andrey S. Trepalin}

\begin{document}

\begin{abstract}
Let $\ka$ be a field of characteristic zero and $G$ be a finite group of automorphisms of projective plane over $\ka$. Castelnuovo's criterion implies that the quotient of projective plane by $G$ is rational if the field $\ka$ is algebraically closed. In this paper we prove that $\Pro^2_{\ka} / G$ is rational for an arbitrary field $\ka$ of characteristic zero.
\end{abstract}

\maketitle
\section{Introduction}

Let $G$ be a finite group and $\ka$ be a field. Consider a pure transcendental extension $\K/\ka$ of transcendence degree $n = \ord G$. Let us identify $\K$ with $\ka\{(x_g)\}$, where $g$ runs through all the elements of the group $G$. The group $G$ naturally acts on $\K$ as $h(x_g)=x_{hg}$. E. Noether \cite{Noe13} asked whether the field of invariants $\K^G$ is rational (i.e. pure transcendental) over $\ka$ or not. In the language of algebraic geometry, this is a question about the rationality of the quotient variety $\A_{\ka}^n/G$.

The most complete answer to this question is known for abelian groups, but even in this case quotient variety can be non-rational. Swan proved that if $G$ is the cyclic group of order $47$ and $\ka=\Q$ then $\K^G$ is not rational (see \cite{Swa69}). A smaller example for the cyclic group of order $8$ was given by Lenstra \cite{Len74}. Further results for abelian groups were obtained in \cite{EM73} and \cite{Vos73}.

For nonabelian groups there are some examples of non-rational field of invariants even in the case $\ka = \kka$. Saltman proved that for any prime $p$ there exists a nonabelian group of order $p^9$ such that $\K^G$ is not rational if $\Char \ka \ne p$ (see \cite{Sal84}). Later this result was developed by Bogomolov, who proved that there exists such a group of order $p^6$ (see \cite{Bog88}), and Moravec, Hoshi and Kang, who proved this result for a group of order $p^5$ (see \cite{Mor12} and \cite{HK11}).

Noether's problem can be generalized as follows. Let $G$ be a finite group, let $V$ be a finite-dimensional vector space over an arbitrary field~$\ka$ and let $\rho: G \rightarrow \mathrm{GL}(V)$ be a representation. The question is if the quotient variety $V/G$ is $\ka$-rational.

Note that $V/G$ has a natural birational structure of a $\A^1$-fibration over $\Pro(V)/G$, which is locally trivial in Zarisky topology. So the rationality of $V/G$ follows from the rationality of $\Pro(V)/G$ (e.g. see \cite[Lemma 1.2]{BK85}).

In this generalization it is natural to start with a low-dimensional case. The most general result is known for dimension $1$ and $2$.

\begin{theorem}[J. L\"uroth]
\label{theorem1}
Let $\ka$ be an arbitrary field and let $G \subset \mathrm{PGL}_2(\ka)$ be a finite subgroup. Then $\Pro^1_{\ka} / G$ is $\ka$-rational.
\end{theorem}

The following theorem is a consequence of the Castelnuovo's rationality criterion.

\begin{theorem}
\label{theorem2}
Let $\ka$ be an algebraically closed field of characteristic zero and let $G \subset \mathrm{PGL}_3(\ka)$ be a finite subgroup. Then $\Pro^2_{\ka} / G$ is $\ka$-rational.
\end{theorem}

On the other hand, if the field $\ka$ is not algebraically closed there is no complete answer. In \cite{Haj87} it is proved that the field $\ka(x,y)^G$ is rational for monomial action of $G$ on the set $\{x, y\}$. It corresponds to $\ka$-rationality of quotients of toric surfaces by groups having an invariant two-dimensional torus on such a surface. From results of the paper \cite{AHK00} it follows that a quotient $\Pro^2_{\ka}/G$ and a quotient $\left( \Pro^1_{\ka} \times \Pro^1_{\ka} \right) / G$ is $\ka$-rational if $G$ is cyclic ($G$ could be infinite).

The main result of this paper is the following.

\begin{theorem}
\label{theorem3}
Let $\ka$ be an arbitrary field of characteristic zero and let $G \subset \mathrm{PGL}_3(\ka)$ be a finite subgroup. Then $\Pro^2_{\ka} / G$ is $\ka$-rational.
\end{theorem}

\begin{corollary}
\label{corollary1}
Let $\ka$ be an arbitrary field of characteristic zero and let $G \subset \mathrm{GL}_3(\ka)$ be a finite subgroup. The field of invariants $\ka(x_1,x_2,x_3)^G$ is $\ka$-rational.
\end{corollary}

For an algebraically closed field $\ka$ it is not known if any quotient of $\Pro^3_{\ka}$ by a finite group is $\ka$-rational (for details see \cite{Pr10}). If in this case the field $\ka$ is not alebraically closed then $\Pro^3_{\ka} / G$ can be non-rational even for abelian group $G$ (see \cite[Example 2.3]{AHK00}).

The plan of proof of \ref{theorem3} is the following. We want to find a normal subgroup $N$ in $G$. If such a group exists then we consider the quotient $\Pro^2_{\ka} / N$. Next, we $G/N$-equivariantly resolve the singularities of $\Pro^2_{\ka} / N$, run the $G/N$-equivariant minimal model program (see \cite{Isk79}) and get a surface $X$. Then we apply the same procedure to the surface $X$ and the group~$G/N$.

This method does not work if the group $G$ is cyclic of prime order or simple nonabelian, but in these cases one can easily prove that the quotient $\Pro^2_{\ka} / G$ is $\ka$-rational (see Subsection \ref{subsection41} and Subsection \ref{subsection45}, respectively).

The structure of the paper is as follows. In Section 2 we describe main notions and some results of the minimal model program which are used in this work. In Section 3 we sketch the classification of finite subgroups in $\mathrm{PGL}_3(\kka)$ where $\kka$ is an algebraically closed field of characteristic zero. In Section 4 we prove Theorem \ref{theorem3}.

The author is grateful to his adviser Yu.\,G.\,Prokhorov for posing the problem and numerous helpful advice and to C.\,A.\,Shramov for useful discussions.

We use the following notation.

\begin{itemize}

\item $\ka$ denotes an arbitrary field of characteristic zero,

\item $\kka$ denotes the algebraic closure of a field $\ka$,

\item $\XX = X \otimes \kka$,

\item $\CG_n$ denotes the cyclic group of order $n$,

\item $\mathfrak{D}_{2n}$ denotes the dihedral group of order $2n$,

\item $\mathfrak{S}_n$ denotes the symmetic group of degree $n$,

\item $\mathfrak{A}_n$ denotes the alternating group of degree $n$,

\item $\omega = e^{\frac{2\pi i}{3}}$,

\item $\operatorname{diag}(\alpha, \beta, \gamma) = \matr{\alpha}{0}{0}{0}{\beta}{0}{0}{0}{\gamma}$,

\item $K_X$ denotes the canonical divisor of a variety $X$,

\item $\Pic(X)$ (resp. $\Pic(X)^G$) denotes the Picard group (resp. $G$-invariant Picard group) of a variety $X$,

\item $\rho(X) = \mathrm{rk} \Pic(X)$, $\rho(X)^G = \mathrm{rk} \Pic(X)^G$,

\item $\F_n$ denotes the rational ruled (Hirzebruch) surface $\Pro_{\Pro^1}(\mathcal{O}~\oplus~\mathcal{O}(n))$,

\item $X \approx Y$ means that $X$ and $Y$ are birationally equivalent.

\end{itemize}

\section{$G$-equivariant minimal model program}

In this section we review main notions and results of $G$-equivariant minimal model program following the papers \cite{Man67}, \cite{Isk79}, \cite{DI1}, \cite{DI2}.

\begin{definition}
\label{definition4}
A {\it rational surface} $X$ is a surface over $\ka$ such that $\XX=X \otimes \kka$ is birationally equivalent to $\Pro^2_{\kka}$.
\end{definition}

\begin{definition}
\label{definition45}
A {\it $\ka$-rational surface} $X$ is a surface over $\ka$ such that $X$ is birationally equivalent to $\Pro^2_{\ka}$.
\end{definition}

\begin{definition}
\label{definition46}
A surface $X$ over $\ka$ is a {\it $\ka$-unirational surface} if there exists a $\ka$-rational variety $Y$ and a dominant rational map $\phi: Y \dashrightarrow X$.
\end{definition}

\begin{definition}
\label{definition5}
A {\it $G$-surface} is a pair $(X, G)$ where $X$ is a projective surface over $\ka$ and $G$ is a subgroup of $\Aut_{\ka}(X)$. A morphism of surfaces $f: X \rightarrow X'$ is called a $G$-morphism $(X, G) \rightarrow (X', G)$ if for each $g \in G$ one has $fg = gf$.
\end{definition}

\begin{definition}
\label{definition6}
A smooth $G$-surface $(X, G)$ is called {\it $G$-minimal} if any birational morphism of smooth $G$-surfaces $(X, G) \rightarrow (X',G)$ is an isomorphism.
\end{definition}

\begin{definition}
\label{definition67}
Let $(X, G)$ be a smooth $G$-surface. A $G$-minimal surface $(Y, G)$ is called a {\it minimal model} of $(X, G)$ if there exists birational $G$-morphism $X \rightarrow Y$.
\end{definition}

The classification of $G$-minimal rational surfaces is well-known due to V.\,Iskovskikh and Yu.\,Manin (see \cite{Isk79} and \cite{Man67}). We introduce some important notions before surveying it.

\begin{definition}
\label{definition7}
A smooth rational $G$-surface $(X, G)$ admits a structure of a {\it conic bundle} if there exists an $G$-morphism $\phi: X \rightarrow C$ such that any scheme fibre is isomorphic to a reduced conic in~$\Pro^2_{\ka}$ and $C$ is a smooth curve.
\end{definition}

\begin{definition}
\label{definition8}
A {\it Del Pezzo surface} is a smooth projective surface $X$ such that the anticanonical divisor $-K_X$ is ample.
\end{definition}

A Del Pezzo surface $\XX$ over $\kka$ is isomorphic to $\Pro^2_{\kka}$, $\Pro^1_{\kka} \times \Pro^1_{\kka}$ or the blowup of $\Pro^2_{\kka}$ at up to 8 points in general position. If $\XX$ is the blowup pf $\Pro^2_{\kka}$ at~$1$ or $2$ points, then $X$ can not be $G$-minimal. The number $K_{\XX}^2$ is called the {\it degree} of a Del Pezzo surface $\XX$.

\begin{theorem}[{\cite[Theorem 1]{Isk79}}]
\label{theorem4}
Let $X$ be a minimal rational $G$-surface. Then either $X$ admits a structure of a conic bundle with $\Pic(X)^{G} \cong \Z^2$, or $X$ is a Del Pezzo surface with $\Pic(X)^{G} \cong \Z$.
\end{theorem}

The following theorem is an important criterion of $\ka$-rationality over an arbitrary perfect field $\ka$.

\begin{theorem}[{\cite[Chapter 4]{Isk96}}]
\label{theorem5}
A minimal rational surface $X$ over a perfect field $\ka$ is $\ka$-rational if and only if the following two conditions are satisfied:

(i) $X(\ka) \ne \varnothing$;

(ii) $K_X^2 \geq 5$.
 \end{theorem}

If a surface $X$ is birationally equivalent to a quotient $\Pro^2_{\ka} / G$ then $X$ is $\ka$-unirational. Thus there exists a Zariski dense subset of $\ka$-points on $X$. The surface $\XX$ is rational by Castelnuovo's criterion. Therefore if $X$ is a smooth surface birationally equivalent to $\Pro^2_{\ka} / G$ and $K_X^2 \geq 5$, then $X$ is $\ka$-rational by Theorem \ref{theorem5}.

An important class of rational surfaces is the class of toric surfaces.

\begin{definition}
\label{definition10}
A \textit{toric variety} is a normal variety over $\ka$ containing an algebraic torus as a Zariski dense subset, such that the action of the torus on itself by left multiplication extends to the whole variety.
\end{definition}

Obviously, a toric variety is $\kka$-rational.

\begin{definition}
\label{definition11}
A variety $X$ is called a \textit{$\ka$-form of a toric variety} if $\XX$ is toric.
\end{definition}

The following lemma is well-known. We give a proof for convenience of the reader.

\begin{lemma}
\label{lemma3}
Let $X$ be a $G$-minimal $\ka$-unirational surface. The following are equivalent:

(i) $X$ is a $\ka$-form of a toric surface;

(ii) $K_X^2 \geq 6$;

(iii) $X$ is isomorphic to $\Pro^2_{\ka}$, a nonsingular quadric $Q \in \Pro^3_{\ka}$, a Del Pezzo surface of degree $6$ or a minimal rational ruled \mbox{surface $\F_n$ ($n \geq 2$)}.
\end{lemma}

\begin{proof}

Let us prove that (i) implies (ii). Assume that $\XX = \Pro^2_{\kka}$, i. e. $X$ is a Brauer--Severi scheme. Since $X(\ka) \ne \varnothing$ one has $X = \Pro^2_{\ka}$ and $K_X^2 = 9$.

From now we assume that $\XX \ne \Pro^2_{\kka}$. Then there exists a birational morphism $f: \XX \rightarrow \left( \F_n \right)_{\kka}$ which is a blowup of $\left( \F_n \right)_{\kka}$ at $k \geq 0$ points $p_1$, \ldots, $p_k$, no two lying in the same fibre of the structure morphism $\pi: \left( \F_n \right)_{\kka} \rightarrow \Pro^1_{\kka}$. Let $\overline{\mathbb{T}}^2 \subset \XX$ be a Zariski dense torus acting on $\XX$. The morphism $f: \XX \rightarrow \left( \F_n \right)_{\kka}$ induces an action of $\overline{\mathbb{T}}^2$ on $\left( \F_n \right)_{\kka}$ and the points \mbox{$p_1$, \ldots, $p_k$} are fixed by this action. Moreover, the morphism \mbox{$\pi f: \XX \rightarrow \Pro^1_{\kka}$} induces an action of $\overline{\mathbb{T}}^2$ on $\Pro^1_{\kka}$ and the points \mbox{$\pi(p_1)$, \ldots, $\pi(p_k)$} are fixed by this action. But a torus acting non-trivially on $\Pro^1_{\kka}$ has exactly two fixed points. Therefore if $\XX$ is toric then $k \leq 2$ and $K_{\XX}^2 = K_X^2 \geq 6$.

Let us prove that (ii) implies (iii). Let $X$ be a $G$-minimal $\ka$-unirational $G$-surface such that $K_X^2 \geq 6$. By Theorem~\ref{theorem4}, either $X$ admits a structure of a conic bundle, or~$X$ is a Del Pezzo surface. If $X$ admits a conic bundle structure and $K_X^2 = 6$ or  $K_X^2 = 7$ then $X$ is not $G$-minimal (see \cite[Theorem~4]{Isk79}). Therefore $K_X^2 = 8$ and $X$ is isomorphic to $\F_n$ ($n \ne 1$) (see \cite[Theorem~3]{Isk79}). A Del Pezzo surface of degree $7$ is never $G$-minimal and a Del Pezzo surface $X$ of degree $8$ is $G$-minimal only if $X$ is isomorphic to a nonsingular quadric $Q \in \Pro^3_{\ka}$ (see \cite[\S3]{Isk79}). Therefore if $X$ is a $G$-minimal Del Pezzo surface and $K_X^2 \geq 6$ then $X$ is isomorphic to $\Pro^2_{\ka}$, a nonsingular quadric $Q \in \Pro^3_{\ka}$ or a Del Pezzo surface of degree $6$.

Let us prove that (iii) implies (i). All surfaces listed in (iii) are $\ka$-forms of a toric surface (see \cite[Section~4, Subsection~6.2]{DI1}).

\end{proof}

The following theorem was proved by Voskresenskii \cite[Theorem~2]{Vos67}. We give another proof in characteristic $0$ based on the theory of $G$-surfaces.

\begin{theorem}
\label{theorem8}
If $X$ is a $\ka$-form of a toric surface such that $X(\ka) \ne \varnothing$ then $X$ is $\ka$-rational.
\end{theorem}

\begin{proof}
If $\widetilde{X} \rightarrow X$ is a minimal resolution of singularities then $\widetilde{X}$ is a $\ka$-form of a toric surface. A minimal model $Y$ of $\widetilde{X}$ is a minimal $\ka$-form of a toric surface, so $K_Y^2 \geq 6$. Therefore $Y$ is $\ka$-rational by Theorem \ref{theorem5}. The surface $Y$ is birationally equivalent to $X$. So $X$ is $\ka$-rational.
\end{proof}

\section{Finite subgroups in $\mathrm{PGL}_3(\kka)$}

In this section we assume that $\ka$ is algebraically closed.

\begin{definition}[{\cite[Chapter I, \S7]{Bl17}}]
\label{definition1}
Any finite subgroup of $\mathrm{GL}_n(\ka)$ is called a {\it linear group} in $n$ variables.
\end{definition}

We will use detailed classification of finite linear subgroups in $3$ variables over $\ka$.

Let a group $G \subset \mathrm{GL}_n(\ka)$ act on the space $V = \ka^n$.

\begin{definition}[{\cite[Chapter I, \S14]{Bl17}}]
\label{definition2}
If the action of the linear group $G$ on $V$ is reducible then the group $G$ is called {\it intransitive}.
Otherwise the group $G$ is called {\it transitive}.
\end{definition}

\begin{definition}[{\cite[Chapter II, \S44]{Bl17}}]
\label{definition3}
Let $G$ be a transitive group. If there exists a non-trivial decomposition $V = V_1 \oplus \dots \oplus V_l$ to subspaces such that for any element $g \in G$ and $i \in \{1, \dots, \mathrm{dim} V\}$ one has \mbox{$gV_i = V_j$} for some $j \in \{1, \dots, \mathrm{dim} V\}$ then the group $G$ is called {\it imprimitive}. Otherwise the group $G$ is called {\it primitive}.
\end{definition}

\begin{convention}
\label{convention2}
Let $f: \mathrm{SL}_n(\ka) \rightarrow \mathrm{PGL}_n(\ka)$ be a natural surjective map. A group $G \subset \mathrm{PGL}_n(\ka)$ is called {\it intransitive } (resp. {\it imprimitive}, {\it primitive}, etc.) if so the group $f^{-1}(G) \subset \mathrm{SL}_n(\ka)$ is.
\end{convention}

The following lemma is well-known.

\begin{lemma}
\label{lemma8}
Any representation of a finite group $G$ in $\mathrm{GL}_n(\ka)$ is conjugate to a representation of the group $G$ in $\mathrm{GL}_n(\overline{\Q}) \subset \mathrm{GL}_n(\ka)$.
\end{lemma}

According to Lemma \ref{lemma8}, it is sufficient to know the classification of finite subgroups of $\mathrm{PGL}_3(\CC)$ to classify finite subgroups of $\mathrm{PGL}_3(\ka)$.

Finite subgroups of $\mathrm{SL}_3(\CC)$ were completely classified in \cite[Chapter V]{Bl17} and \cite[Chapter XII]{MBD16}. Using this result and Lemma \ref{lemma8} one can obtain the following theorem.

\begin{theorem}
\label{theorem7}
Any finite subgroup in $\mathrm{PGL}_3(\ka)$ modulo conjugation is one of the following.

\medskip

Intransitive groups:

\medskip

{\rm{(A)}} A diagonal abelian group.

\smallskip

{\rm{(B)}} A group having a unique fixed point on $\Pro^2_{\ka}$.

\medskip

Imprimitive groups:

\medskip

{\rm{(C)}} A group having a normal diagonal abelian subgroup $N$ such that $G / N \cong \CG_3$.

\smallskip

{\rm{(D)}} A group having a normal diagonal abelian subgroup $N$ such that $G / N \cong \mathfrak{S}_3$.

\medskip

Primitive groups having a non-trivial normal subgroup (the Hessian group and its subgroups):

\medskip

{\rm{(E)}} The group $\CG_3^2 \rtimes \CG_4$ of order $36$.

\smallskip

{\rm{(F)}} The group $\CG_3^2 \rtimes Q_8$ of order $72$.

\smallskip

{\rm{(G)}} The Hessian group $\CG_3^2 \rtimes \mathrm{SL}_2(\F_3)$ of order $216$.

\medskip

Simple groups:

\medskip

{\rm{(H)}} The icosahedral group $\mathfrak{A}_5$ of order $60$.

\smallskip

{\rm{(I)}} The Klein group $\mathrm{PSL}_2(\F_7)$ of order $168$.

\smallskip

{\rm{(K)}} The Valentiner group $\mathfrak{A}_6$ of order $360$.

\end{theorem}

\section{Rationality of the quotient variety}

In this section, for each finite group $G \subset \mathrm{PGL}_3(\ka)$ (see Theorem~\ref{theorem7}), we prove that the quotient variety $\Pro^2_{\ka} / G$ is $\ka$-rational.

We will use the following definition for convenience.

\begin{definition}
\label{definition9}
Let $S$ be a $G$-surface, $\widetilde{S} \rightarrow S$ be its ($G$-equivariant) minimal resolution of singularities, and $Y$ be a $G$-equivariant minimal model of $\widetilde{S}$. We call the surface $Y$ a \textit{$G$-MMP-reduction} of $S$.
\end{definition}

The $\ka$-rationality of $\Pro^2_{\ka} / G$ is proved in Subsection \ref{subsection41} and Subsection \ref{subsection45} in cases where $G$ is conjugate to a diagonal abelian or a simple nonabelian group respectively over $\kka$.

Otherwise there is a non-trivial normal subgroup $N \lhd G$. Then we use the following lemma.

\begin{lemma}
\label{lemma1}
Let $X$ be a $G$-surface, $N$ be a normal subgroup of $G$ and $Y$ be a $(G / N)$-MMP-reduction of $X / N$. The surfaces $X / G$ and $Y / ( G / N)$ are birationally equivalent over $\ka$. 
\end{lemma}

\begin{proof}
One has $X / G = (X / N) / (G / N)$. The birational map \mbox{$X / N \dashrightarrow Y$} is $(G / N)$-equivariant. Thus $(X / N) / (G / N)$ and $Y / (G / N)$ are birationally equivalent over $\ka$.
\end{proof}

Therefore if $X$ is a $(G / N)$-MMP-reduction of $\Pro^2_{\ka} / N$ then $\ka$-rationality of $X / (G / N)$ implies $\ka$-rationality of $\Pro^2_{\ka} / G$. If there is a normal subgroup $M$ in the group $G / N$ we can repeat this procedure applying Lemma \ref{lemma1} to the surface $X$ and the normal subgroup $M \lhd G / N$.

\subsection{Diagonal abelian groups}
\label{subsection41}
Let $G \subset \mathrm{PGL}_3(\ka)$ be an abelian subgroup which is conjugate to a diagonal subgroup of $\mathrm{PGL}_3(\kka)$. Then the action of $G$ on $\Pro^2_{\kka}$ can be considered as an action of a finite subgroup of a torus in $\Pro^2_{\kka}$.

The following lemma is well-known.

\begin{lemma}
\label{lemma9}
Let $\XX$ be an $n$-dimensional toric variety over a field $\kka$ of arbitrary characteristic and let $G$ be a finite subgroup conjugate to a subgroup of $n$-dimensional torus $\overline{\mathbb{T}}^n \subset \XX$ acting on $\XX$. Then the quotient $\XX / G$ is a toric variety.

In particular, if $G$ is a finite cyclic subgroup of connected component of the identity $\Aut^0(\XX) \subset \Aut(\XX)$, then the quotient $\XX / G$ is a toric variety.
\end{lemma}

\begin{proof}
The algebra of regular functions of $\overline{\mathbb{T}}^n$ is
$$
\kka[\overline{\mathbb{T}}^n] = \kka \left[ x_1, \dots, x_n, \frac{1}{x_1}, \dots, \frac{1}{x_n} \right]
$$
and its monomials form a lattice $M \cong \Z^n$. The group $G$ acts on the algebra $\kka[\overline{\mathbb{T}}^n]$ and if $P \in \kka[\overline{\mathbb{T}}^n]$ is $G$-invariant polynomial then~$P$ consists of $G$-invariant monomials. Therefore $G$-invariants form an \mbox{$n$-dimensional} sublattice in the lattice $M$. It means that $\overline{\mathbb{T}}^n / G$ is an $n$-dimensional Zariski dense torus in $\XX / G$. For any $\overline{\mathbb{T}}^n$-orbit $\overline{A} \subset \XX$ the torus $\overline{\mathbb{T}}^n/ G$ acts on $\overline{A} / G$. So $\XX / G$ is a toric variety.

If $G$ is a finite cyclic subgroup of connected component of the identity in $\Aut(\XX)$ then $G$ is a subgroup of a maximal torus in $\Aut(\XX)$ (see \cite[Theorem 11.10]{Bor91}). Therefore $\XX / G$ is a toric variety.
\end{proof}

By Lemma \ref{lemma9} the surface $\Pro^2_{\ka} / G$ is a $\ka$-form of a toric surface. Thus it is $\ka$-rational by Theorem \ref{theorem8}.

\subsection{Groups having a unique fixed point}

The fixed $\kka$-point $p$ of $G$ is defined over $\ka$ because it is unique.
Let $V \rightarrow \Pro^2_{\ka}$ be the blowup at $p$. The surface $V$ admits a $G$-equivariant $\Pro^1_{\ka}$-bundle structure $V \rightarrow \Pro^1_{\ka}$ whose fibres are proper transforms of lines passing through the point $p$. Obviously this $\Pro^1_{\ka}$-bundle has a $G$-invariant section: the exceptional divisor of the blowup at $p$. So one has (e.g. see \cite[Lemma 1.2]{BK85})
$$
\Pro^2_{\ka} / G \approx V / G \approx \A^1_{\ka} \times \Pro^1_{\ka} / G
$$
and it is $\ka$-rational by Theorem \ref{theorem1}.

\subsection{Imprimitive groups}
Every imprimitive group $G \subset \mathrm{PGL}_3(\ka)$ contains a normal abelian subgroup $N$ which is conjugate to a diagonal subgroup in $\mathrm{PGL}_3(\kka)$. The quotient group $G/N$ permutes eigenspaces of $N$ transitively so it is isomorphic to $\CG_3$ or~$\mathfrak{S}_3$ (it corresponds to cases {\rm{(C)}} and {\rm{(D)}} of Theorem \ref{theorem7}). Moreover, a $G / N$-MMP-reduction of $\Pro^2_{\ka} / N$ is a $\ka$-form of a toric surface $X$ by Lemma~\ref{lemma9}. Therefore the $\ka$-rationality of $\Pro^2_{\ka} / G$ follows from Lemma \ref{lemma1} and the following proposition:

\begin{proposition}
\label{proposition1}
Let $X$ be a $G$-minimal $\ka$-unirational $\ka$-form of a toric surface and let $G$ be a group $\CG_3$ or $\mathfrak{S}_3$ acting on $X$. Then $X / G$ is $\ka$-rational.
\end{proposition}

In the proof of this proposition quotient singularities of certain types play an important role. The following remark is useful.

\begin{remark}
\label{remark2}
Let $\xi$ be a primitive $n$-th root of unity and let the group $\CG_n$ act on a smooth surface $X$ with a fixed point $P \in X$. Let $f: X \rightarrow X/ \CG_n$ be the quotient map.
If $\CG_n$ acts on the tangent space at the point $P$ as $\operatorname{diag}(\xi_n, \xi_n^{n-1})$ then the point $f(P)$ is a Du Val singularity of type~$A_{n-1}$.

If a surface $S$ has only Du Val singularities and $\pi: \widetilde{S} \rightarrow S$ is a minimal resolution of singularities then $K_{\widetilde{S}} = \pi^* K_S$.
\end{remark}

\begin{proposition}
\label{proposition2}
Let a group $G$ contain a normal subroup $\CG_p$, where $p$ is prime. If $X$ is a $G$-minimal $\ka$-unirational $\ka$-form of a toric surface then there exists a $G / \CG_p$-MMP-reduction $Y$ of $X / \CG_p$ such that $Y$ is a $\ka$-form of a toric surface. In particular, $X / \CG_p$ is $\ka$-rational.
\end{proposition}

\begin{proof}

By Lemma \ref{lemma3} the surface $X$ is isomorphic to $\Pro^2_{\ka}$, a quadric in $\Pro^3_{\ka}$, $\F_n$ ($n \geq 2$) or a Del Pezzo surface of degree $6$.
If the group $\CG_p$ is a subgroup of the connected component of the identity $\Aut^0(\XX) \subset \Aut(\XX)$ then any $G / \CG_p$-MMP-reduction $Y$ of $X / \CG_p$ is a $\ka$-form of a toric surface by Lemma \ref{lemma9}.
Otherwise $X$ is isomorphic to a quadric in $\Pro^3_{\ka}$ or a Del Pezzo surface of degree $6$ and the group $\CG_p$ acts on $\Pic(X)$ faithfully.

If $\XX = \Pro^1_{\kka} \times \Pro^1_{\kka}$ and $\CG_p$ is not a subgroup of $\Aut^0(\XX)$ then $p = 2$ and the action of the group $\CG_2$ is conjugate to
$$
(x_1 : x_0 ; y_1 : y_0) \mapsto (y_1 : y_0 ; x_1 : x_0)
$$
where $(x_1 : x_0 ; y_1 : y_0)$ are homogeneous coordinates on $\Pro^1_{\kka} \times \Pro^1_{\kka}$. The set of fixed points is a curve $\frac{x_1}{x_0}=\frac{y_1}{y_0}$ on $\Pro^1_{\kka} \times \Pro^1_{\kka}$ whose class in \mbox{$\Pic(\Pro^1_{\kka} \times \Pro^1_{\kka})$} equals to $-\frac{1}{2} K_{\Pro^1_{\kka} \times \Pro^1_{\kka}}$. The quotient variety $\overline{Y} = (\Pro^1_{\kka} \times \Pro^1_{\kka}) / \CG_2$ is nonsingular. By the Hurwitz formula one has $K_{\overline{Y}}^2=\frac{1}{2} (\frac{3}{2}K_{\XX})^2 = 9$. The surface $\overline{Y}$ is rational so it is isomorphic to $\Pro^2_{\kka}$. Thus $Y$ is a $\ka$-form of a toric surface.

If $X$ is a Del Pezzo surface of degree $6$ then $\XX$ is isomorphic to a blowup of $\Pro^2_{\kka}$ at three points in general position. The automorphism group of $\XX$ is $\overline{\mathbb{T}}^2 \rtimes \mathfrak{D}_{12}$ where $\overline{\mathbb{T}}^2$ is a two-dimensional torus over $\kka$ and $\mathfrak{D}_{12}$ faithfully acts on $\Pic(\XX)$ and the set of $(-1)$-curves $\{ C_i \}$, \mbox{$1 \leq i \leq 6$} forming a hexagon (the irreducible components of the exceptional divisor of the blowup and the proper transforms of lines passing through a pair of points of the blowup).

Recall that $\mathfrak{D}_{12}$ generated by two elements $r$ and $s$ with relations
$$
r^6=s^2=(rs)^2=1.
$$
There are four cyclic subgroups of prime order in the group $\mathfrak{D}_{12}$ up to conjugacy:

(a) $\CG_2 = <s>$;

(b) $\CG_2 = <rs>$;

(c) $\CG_2 = <r^3>$;

(d) $\CG_3 = <r^2>$.

\begin{center}
\includegraphics[scale=0.8]{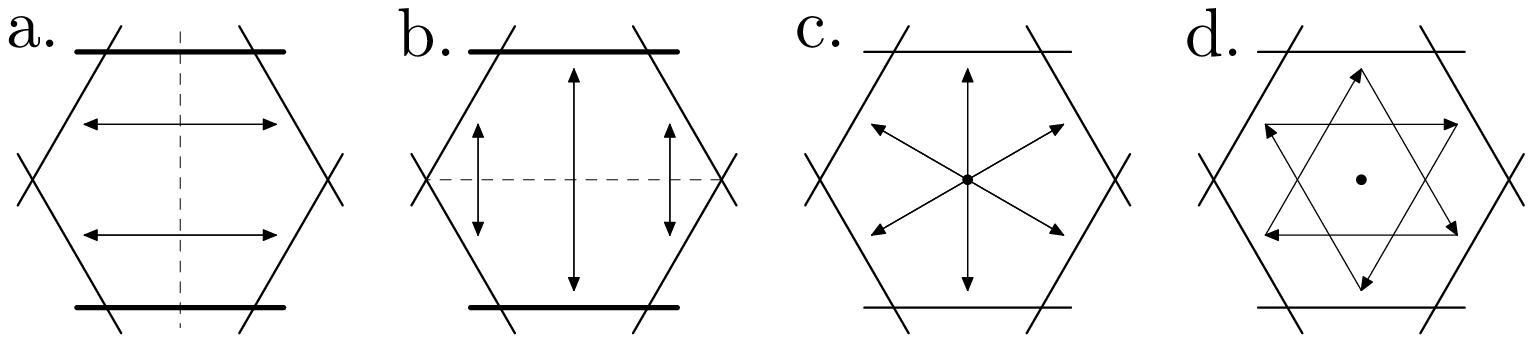}
\end{center}

In the case (a) we have two $\CG_2$-invariant disjoint $(-1)$-curves and the others are not $\CG_2$-invariant. Therefore the action of $G$ is not minimal.

In the case (b) we have the $\CG_2$-invariant pair of disjoint $(-1)$-curves (two other $\CG_2$-invariant pairs of $(-1)$-curves are not disjoint). As above the action of $G$ is not minimal.

In the remaining cases (c) and (d) let us denote the quotient morphism $X \rightarrow X / \CG_p$ by $g$. The group $\CG_p$ has no fixed points on $\sum \limits_{i=1}^6 C_i$. Since $\sum \limits_{i=1}^6 C_i \sim -K_X$ is ample, $\CG_p$ has no curves of fixed points. Hence $K_X = g^*K_{X / \CG_p}$ and so $-K_{X / \CG_p}$ is ample.

In the case (c) by the Lefschetz fixed-point formula there are four fixed points. So the quotient $\XX / \CG_2$ has exactly four $A_1$-singularities, $K_{X / \CG_2}^2 = 3$ and $X / \CG_2$ is isomorphic to a singular cubic surface in $\Pro^3_{\ka}$. Such a cubic is $G / \CG_2$-equivariantly birationally equivalent to $\Pro^2_{\ka}$ (see \cite[Lemma 1.1]{CT88}).

In the case (d) a generator $h \in \CG_3$ acts on the tangent space of a fixed point in one of the following ways: $\operatorname{diag}(\omega, \omega)$, $\operatorname{diag}(\omega^2, \omega^2)$ or $\operatorname{diag}(\omega, \omega^2)$. Denote by $a$, $b$ and $c$ the number of corresponding points. Then the Lefschetz fixed-point formula and the holomorphic Lefschetz fixed-point formula yield
$$
a + b + c = 3, \quad -\frac{a}{3 \omega} -\frac{b}{3 \omega^2} + \frac{c}{3} = 1.
$$
Therefore $a = b = 0$, $c = 3$ and so $\XX / \CG_3$ has exactly three \mbox{$A_2$-singularities}. One has $K_{X / \CG_3}^2 = 2$. Thus $X / \CG_3$ is a singular Del Pezzo surface of degree $2$ with Du Val singularities. The image $g \left( \sum \limits_{i=1}^6 C_i \right)$ consist of two irreducible components which are two $(-1)$-curves contained into the smooth locus of $\XX / \CG_3$. Moreover, $\rho(\XX / \CG_3) = \rho(\XX)^{\CG_3} = 2$. Therefore contractions of these curves correspond to two extremal rays of the Mori cone of $\XX / \CG_3$ so all other curves on $\XX / \CG_3$ have positive selfintersection.

For any Del Pezzo surface $V$ of degree $2$ with at worst Du Val singularities the linear system $|-K_V|$ is base point free and defines a double cover $f: V \rightarrow \Pro^2_{\kka}$ branched over a reduced quartic $B \subset \Pro^2_{\kka}$. In our case $V = \XX / \CG_3$ from the local equations one can obtain that $B$ is irreducible and has $3$ cusps. Let us consider three lines $l_1$, $l_2$, $l_3$ on~$\Pro^2_{\kka}$ passing through pairs of these cusps. The preimage $f^{-1}(l_i)$ is a pair of rational curves meeting each other at two singular points. Let $a$ be one of the irreducible components of $f^{-1}(l_i)$. Then one has:
$$
a \cdot K_{\XX / \CG_3} = \frac{1}{2} \cdot \left(a + i(a) \right) \cdot K_{\XX / \CG_3} = \frac{1}{2} \left( -K_{\XX / \CG_3}^2 \right) = -1,
$$
$$
a \cdot f^{-1}(l_i) = a \cdot \left( -K_{\XX / \CG_3} \right) = 1,
$$
$$
a \cdot i(a) < \frac{1}{2} \left( a + i(a) \right)^2 = \frac{1}{2} \left( -K_{\XX / \CG_3} \right)^2 = 1.
$$
Therefore the intersection number of each pair of distinct irreducible components of $f^{-1}(l_1)$, $f^{-1}(l_2)$ and $f^{-1}(l_3)$ is positive and less than~$1$. If $\pi: \widetilde{X / \CG_3} \rightarrow X / \CG_3$ is the minimal $G / \CG_3$-equivariant resolution of singularities then the proper transform $\pi^{-1}_* f^{-1} (l_1 + l_2 + l_3)$ consist of six disjoint \mbox{$(-1)$-curves}, defined over~$\ka$. We can $G / \CG_3$-equivariantly contract these six curves and get a surface $Y$ with
$$
K_Y^2 = K_{\widetilde{X / \CG_3}}^2 + 6 = K_{X / \CG_3}^2 + 6 = \frac{1}{3} K_X^2 + 6 = 8.
$$
So $Y$ is a $\ka$-form of a toric surface by Lemma \ref{lemma3}.

\end{proof}

Now we come to the proof of Proposition \ref{proposition1}.

\begin{proof}[Proof of Proposition \ref{proposition1}]

If $G = \CG_3$ then Proposition \ref{proposition1} directly follows from Proposition \ref{proposition2}.

If $G = \mathfrak{S}_3$ then we apply Lemma \ref{lemma1} to the surface $X$ and the normal subgroup $\CG_3 \lhd \mathfrak{S}_3$. By Proposition \ref{proposition2} there exists a \mbox{$\CG_2$-MMP-reduction}~$Y$ of $X / \CG_3$ which is a $\ka$-form of a toric surface. The surface $Y / \CG_2$ is $\ka$-rational by Proposition \ref{proposition2}.

\end{proof}

\begin{corollary}
\label{corollary2}
Let $X$ be a Del Pezzo surface of degree $6$ over $\ka$, $X(\ka) \ne \varnothing$ and $G$ be a finite subgroup of automorphisms of $X$. The quotient variety $X / G$ is $\ka$-rational.
\end{corollary}

\begin{proof}
The proof is a combination of Lemma \ref{lemma9}, Proposition \ref{proposition2} and Lemma \ref{lemma1} applied to the normal subgroups $N \lhd G$, $M \lhd G / N$, \mbox{$L \lhd (G / N) / M$} and $K = ((G / N) / M) / L$ defined in the following way:

$$
N = G \cap \overline{\mathbb{T}}^2 \lhd G \subset \overline{\mathbb{T}}^2 \rtimes \mathfrak{D}_{12},
$$
$$M = G / N \cap \CG_2 \lhd G / N \subset \mathfrak{D}_{12},
$$
$$
L = (G / N) / M \cap \CG_3 \lhd (G / N) / M \subset \mathfrak{S}_3 \cong \mathfrak{D}_{12} / \CG_2,
$$
$$
K = ( (G / N) / M) / L \lhd \CG_2 \cong \mathfrak{S}_3 / \CG_3.
$$
\end{proof}

\subsection{Primitive groups having a non-trivial normal subgroup}

Primitive groups having non-trivial normal subgroups are the Hessian group ${\rm{(G)}} = \CG_3^2 \rtimes \mathrm{SL}_2(\F_3)$ and its subgroups ${\rm{(F)}} = \CG_3^2 \rtimes Q_8$ and \mbox{${\rm{(E)}} = \CG_3^2 \rtimes \CG_4$} (for details see \cite[Section 4]{AD06}). Moreover, we may assume that they have a common normal subgroup $N$ isomorphic to~$\CG_3^2 \rtimes \CG_2$. The quotients of the groups ${\rm{(E)}}$, ${\rm{(F)}}$, ${\rm{(G)}}$ by $N$ are as follows: ${\rm{(E)}}/N = \CG_2$, ${\rm{(F)}}/N = \CG_2^2$, ${\rm{(G)}}/N = \mathfrak{A}_4$.

The group $N$ is generated by elements of order $2$. Any element $g$ of order $2$ in $\mathrm{PGL}_3(\kka)$ is conjugate to $\operatorname{diag}(-1, 1, 1)$, that is, $g$ is a reflection. The quotient $\Pro^2_{\kka} / N$ is a weighted projective space by the Chevalley--Shephard--Todd theorem (see \cite{ST54}). Any weighted projective space is a toric variety, therefore any $G/N$-MMP-reduction $X$ of $\Pro^2_{\ka} / N$ is a $\ka$-form of a toric surface. By Lemma \ref{lemma1} the $\ka$-rationality of the quotients $\Pro_{\ka}^2 / {\rm{(E)}}$ and $\Pro_{\ka}^2 / {\rm{(F)}}$ follows from the $\ka$-rationality of the quotients $X / \CG_2$ and $X / \CG_2^2$.

The $\ka$-rationality of $X / \CG_2$ directly follows from Proposition \ref{proposition2}.

The $\ka$-rationality of $X / \CG_2^2$ follows from Lemma \ref{lemma1} applied to the surface $X$ and a normal subgroup $\CG_2 \lhd \CG_2^2$. By Proposition \ref{proposition2} there exists a $\CG_2$-MMP-reduction $Y$ of $X / \CG_2$ which is a $\ka$-form of a toric surface. The surface $Y / \CG_2$ is $\ka$-rational by Proposition \ref{proposition2}.

The Hessian group ${\rm{(G)}}$ is generated by reflections of order $3$ (see e.g. \cite[Section 4]{AD06}). Therefore the quotient $\Pro^2_{\kka} / {\rm{(G)}}$ is a weighted projective space by the Chevalley--Shephard--Todd theorem. It is a toric surface so $\Pro^2_{\ka} / {\rm{(G)}}$ is $\ka$-rational by Theorem \ref{theorem8}.

\subsection{Simple groups}
\label{subsection45}
Let $G \subset \mathrm{PGL}_3(\ka)$ be a simple group. Note that each simple group {\rm{(H)}}--{\rm{(K)}} from Theorem \ref{theorem7} is generated by elements of order $2$ because its order must be even (so there is at least one element of order $2$) and elements of order $2$ generate a normal subgroup. Therefore, similar to the previous subsection, $G$ is generated by reflections and the quotient $\Pro^2_{\kka} / G$ is a weighted projective space. It is a toric surface so $\Pro^2_{\ka} / G$ is $\ka$-rational by Theorem \ref{theorem8}.

\bibliographystyle{alpha}
\bibliography{my_ref}
\end{document}